\documentclass[12pt]{article}
\newtheorem{theorem}{Theorem}
\newtheorem{lemma}{Lemma}
\newtheorem{conj}{Conjecture}
\newtheorem{cor}{Corollary}
\usepackage{amssymb}
\begin{document}
\title{Colouring planar graphs with a precoloured induced cycle}
\author {Ajit A. Diwan \\ Department of Computer Science and Engineering,\\
Indian Institute of Technology Bombay, Mumbai 400076, India.\\
email: \texttt{aad@cse.iitb.ac.in}}
\maketitle
\abstract{ Let $C$ be a cycle and $f : V(C) \rightarrow \{c_1,c_2,\ldots,c_k\}$ a proper 
$k$-colouring of $C$ for some $k \ge 4$. We say the colouring $f$ is \emph{safe} if for 
any planar graph $G$ in which $C$ is an induced cycle, there exists a proper $k$-colouring 
$f'$ of $G$ such that $f'(v) = f(v)$ for all $v \in V(C)$. The only safe $4$-colouring is any 
proper colouring of a triangle. We give a simple necessary condition for a $k$-colouring 
of a cycle to be safe and conjecture that it is sufficient for all $k \ge 4$. The sufficiency
for $k=4$ follows from the four colour theorem and we prove it for $k = 5$, independent of the
four colour theorem. We show that a stronger condition is sufficient for all $k \ge 4$. As 
a consequence, it follows that any proper $k$-colouring of a cycle that uses at most $k-3$ distinct 
colours is safe. Also, any proper $k$-colouring of a cycle of length at most $2k-5$ that uses at 
most $k-1$ distinct colours is safe.}

\section{Introduction}

A proper $k$-colouring of a graph $G$ is an assignment $f : V(G) \rightarrow \{c_1,c_2,\ldots,c_k\}$
of colours to the vertices of $G$ such that for every edge $uv$ in the graph $f(u) \neq
f(v)$. We will refer to proper colourings as simply colourings and a graph is $k$-colourable if
there exists a $k$-colouring. The study of colourings of planar graphs has a long history, starting 
with the four colour conjecture and its eventual proof~\cite{AH1}~\cite{AH2}~\cite{RSST}. Many 
variations of this have since been considered, the most important of which is perhaps the notion of 
list-colouring. Let $L(v)$ be a list of allowed colours assigned to each vertex $v \in V(G)$. 
An $L$-colouring of $G$ is a proper colouring such that each vertex $v$ is assigned a colour in $L(v)$. 
The usual colouring problem is the special case when the lists $L(v)$ are the same for all vertices. 
While all planar graphs are four colourable, it is not true that they are $L$-colourable if $|L(v)| 
\ge 4$ for all vertices $v$~\cite{V}. Thomassen~\cite{T} showed that all planar graphs are $L$-colourable 
if $|L(v)| \ge 5$ for each vertex $v$. Subsequently, there has been a lot of work on finding conditions 
on the lists $L(v)$ that guarantee any planar graph is $L$-colourable. A recent result of Zhu~\cite{Z}
shows that if $|L(v)|=4$ for each vertex $v$, and $|L(u) \cap L(v)| \le 2$ for any edge $uv$,
then any planar graph is $L$-colourable. Many such results with different restrictions were
proved earlier, and several conjectured, as mentioned in ~\cite{Z}. 

Another variation considered is colouring or $L$-colouring with some precoloured vertices,
whose colours are specified. This is equivalent to considering an $L$-colouring in which some 
vertices have lists of size 1. This can also be viewed as extending a partial assignment
of colours to some vertices to a colouring or $L$-colouring of the whole graph.  
Albertson~\cite{A} showed that every planar graph with precoloured vertices is five colourable 
if the precoloured vertices are at distance at least 4 from each other and the precolouring uses 
at most five distinct colours. No such result is possible for four colours, even if only two vertices 
are precoloured. The corresponding result for $L$-colouring  was shown in ~\cite{DLMP}. There exists 
a constant $M$ such that if the precoloured vertices are at distance at least $M$ from each other, 
and $|L(v)|=5$ for all other vertices, then any planar graph is $L$-colourable. 

In most of these results, the restrictions on the lists are in terms of their sizes, or
the sizes of their intersections, but not the actual elements in the list. Here, we consider
another special case where the vertices of a connected induced subgraph are assigned lists of size 1 
and all other vertices have the same list of a specified size. This is equivalent to asking when can
a specified $k$-colouring of a graph $H$ be extended to a $k$-colouring of any planar graph
that contains $H$ as in induced subgraph.  We call such a $k$-colouring of the graph $H$
a \emph{safe} $k$-colouring. The subgraph we consider is a cycle. If $H$ is any $2$-connected planar
graph, a $k$-colouring of $H$ is safe if and only if the colouring of every non-separating induced
cycle in $H$ is safe. Thus considering a cycle as the subgraph $H$ is a natural choice.

We give a simple necessary condition for a $k$-colouring of a cycle to be safe. We conjecture
that this condition is sufficient for all $k \ge 4$. The sufficiency for $k=4$ follows from the 
four colour theorem. We prove it for $k=5$. The proof is independent of and much simpler
than the proof of the four colour theorem. We prove that a stronger condition is sufficient
for all $k \ge 4$. As a consequence, any $k$-colouring of a cycle that uses at most
$k-3$ distinct colours is safe. Also, any $k$-colouring of a cycle of length at most $2k-5$
that uses at most $k-1$ colours is safe. If the conjecture is true, then any $k$-colouring of
a cycle of length at most $3k-10$ that uses at most $k-2$ colours is safe. However, we have
not been able to prove this. 

Such results have been considered for 3-colourings of planar triangle-free graphs.
In ~\cite{DL} for example, 3-colourings of a cycle of length at most 8 contained in a planar
triangle-free graph that cannot be extended to a 3-colouring of the whole graph have been
characterized. The same has been done for cycles of length 9 in ~\cite{CEHL}. However, we do not 
know any such results for general planar graphs when the number of colours and the cycle length is 
arbitrary. Also, these results are based on characterizing all critical graphs, that is graphs
for which the colouring of the cycle cannot be extended to the whole graph, but can be extended
for any proper subgraph that contains the cycle. We do not attempt to do this here. To show
that a colouring is not safe, it is enough to show one graph for which the colouring of the
cycle cannot be extended.  To show it is safe, we use properties of the colouring to extend it
to the whole graph.

\section{Main Result}

Let $C = v_1,v_2,\ldots,v_l$ be a cycle of length $l \ge 3$ and $f$ a safe $k$-colouring of $C$. 
Any $k$-colouring of $C$ obtained from $f$ by permuting the colours is also safe. The 
colouring $f_i$ defined by $f_i(v_j) = f(v_{i+j})$, where addition is modulo $l$, is also safe, 
for all $1 \le i < l$. Similarly, the colouring obtained by reversing the cycle, that is $f_r(v_i) = 
f(v_{l+1-i})$ is safe.  We say these colourings of $C$ are equivalent and consider safety to be 
a property of the equivalence class. 

Let $[i,j] = \{i, i+1, \ldots, j\}$ for $1 \le i,j \le l$, where again addition is considered modulo 
$l$. Let $[i,j) = [i,j] \setminus \{j\}$, $(i,j] = [i,j] \setminus \{i\}$, and $(i,j) = [i,j] \setminus
\{i,j\}$. Let $F[i,j] = \{ c_t \mid 1 \le t \le k, \exists  m \in [i,j] \ f(v_m) = c_t\}$ be the subset 
of colours that occur in the subpath of $C$ from $v_i$ to $v_j$. The sets $F[i,j), F(i,j], F(i,j)$ are 
defined similarly.

To prove that a $k$-colouring $f$ of $C$ is safe, we consider an arbitrary planar graph $G$ such that 
$C$ is an induced cycle in $G$, and show that the colouring of $C$ can be extended to a 
$k$-colouring of $G$. It is sufficient to do this in the case $C$ is a non-separating induced cycle 
in $G$, since otherwise we can consider each component of $G-V(C)$ separately. Thus we may assume 
$C$ bounds a face of $G$ and without loss of generality, we can embed $G$ so that $C$ is the boundary 
of the external face. We may  further assume $G-V(C)$ is not empty, and add edges to $G$ that are 
not chords of $C$ if needed, so that every internal face of $G$ is bounded by a triangle. We call 
such a graph $G$ a \emph{chordless near-triangulation}. Thus a $k$-colouring $f$ of a cycle $C$ is safe 
if and only if for any chordless near-triangulation $G$ with external face bounded by the cycle $C$, 
there exists a $k$-colouring of $G$ that extends $f$.

We first give some simple necessary conditions for a $k$-colouring to be safe.

\begin{lemma}
\label{unsafe}
Let $f$ be a $k$-colouring of a cycle $C$ of length $l$. If $f$ satisfies any one
of the following conditions, then $f$ is not safe.
\begin{enumerate}
\item
$|F[1,l]| = k$.
\item
There exist indices $1 \le p < q \le l$ such that $|F[p,q] \cap F[q,p]| \ge k-1$.
\item
There exist indices $1 \le p < q < r \le l$ such that $|F[p,q] \cap F[q,r] \cap F[r,p]| \ge k-2$.
\end{enumerate}
\end{lemma}

\noindent {\bf Proof:} If $|F[1,l]| = k$ then the $k$-colouring cannot be extended to a $k$-colouring  
of the wheel $W_l$, hence it is not safe. If there are indices $p < q$ such that $|F[p,q] \cap F[q,p]| 
\ge k-1$, construct a chordless near-triangulation by adding a vertex $u$ that is adjacent to all 
vertices $v_m$ for $m \in [p,q]$, and a vertex $v$ that is adjacent to $u$ and all vertices $v_m$ for
$m \in [q,p]$. The $k$-colouring of the cycle cannot be extended to this near-triangulation, hence such a 
colouring is not safe. The same argument holds if there are indices $p < q < r$ such that $|F[p,q] \cap 
F[q,r] \cap F[r,p]| \ge k-2$. Construct a near-triangulation by adding a triangle $u,v,w$ with $u$ adjacent 
to vertices $v_m$ for $m \in [p,q]$, $v$ adjacent to vertices $v_m$ for $m \in [q,r]$ and $w$ adjacent to 
vertices $v_m$ for $m \in [r,p]$. Again it is easy to check that the colouring of the cycle cannot be extended 
to the near-triangulation, hence any such $k$-colouring is not safe. 
\hfill $\Box$

We call a $k$-colouring of a cycle that satisfies any of the conditions in Lemma~\ref{unsafe}
a \emph{bad} $k$-colouring. If $k=4$ and $l \ge 4$, then every  $4$-colouring of $C$ is bad. Consider 
any $4$-colouring of a cycle $C$ of length at least $4$. Assume $|F[1,l]| < 4$ otherwise condition 1 
in Lemma~\ref{unsafe} is satisfied. If $|F[1,l]| = 2$, we have $|F[1,2]| = 2$, and $F[1,2] = F[2,3] = 
F[3,1]$, which implies $f$ satisfies condition 3 in Lemma~\ref{unsafe} with $p = 1, q = 2, r = 3$. 
If $|F[1,l]| = 3$, there exists an index $i$ such that $|F[i,i+2]| = 3$. If $|F[i+2,i]| = 3$, then $f$ 
satisfies condition 2 in Lemma~\ref{unsafe} with $\{p,q\} = \{i,i+2\}$, otherwise we have $l \ge 5$, 
$f(v_{i+3}) = f(v_i)$ and $f(v_{i+4}) = f(v_{i+2})$. In this case, we have $|F[i+1,i+3]| = |F[i+3,i+1]| = 3$ 
and $f$ again satisfies condition 2 in Lemma~\ref{unsafe} with $\{p,q\} =  \{i+1,i+3\}$. Therefore, 
the only $4$-colourings that are not bad are those of the triangle. It follows from the four colour 
theorem that any $4$-colouring of a triangle is safe and hence a $4$-colouring is safe if and only if 
it is not bad. We conjecture that this property holds for all $k \ge 4$.

\begin{conj}
\label{safe}
A $k$-colouring of a cycle is safe if and only if it is not bad, for all $k \ge 4$.
\end{conj}

We prove Conjecture~\ref{safe} for $k = 5$. Although the statement of Conjecture~\ref{safe} may be
considered to be a generalization of the four colour theorem, the proof for $k=5$ is much
simpler and does not depend on the four colour theorem. We in fact show that $k$-colourings
satisfying a stronger property are safe for all $k \ge 4$. All 
$5$-colourings that are not bad satisfy the stronger property and are therefore safe. Moreover,
there exist such $k$-colourings satisfying the stronger property for cycles of all lengths, for
all $k \ge 5$. However, this proof depends crucially on the fact that $k \ge 5$, and cannot be 
easily adapted to prove the four colour theorem itself. Also, for $k \ge 6$, there  are $k$-colourings 
that are not bad but do not satisfy the stronger property, so this does not prove Conjecture~\ref{safe}
for $k \ge 6$. 

Henceforth, we will assume that $k \ge 5$ is a fixed integer. We call a $k$-colouring
of a cycle of length $l$ a \emph{good} $k$-colouring if it satisfies one of the following properties.
\begin{enumerate}
\item
$|F[1,l]| \le k-3$.
\item
$|F[1,l]| = k-2$ and there exist indices $1 \le p < q \le l$ and a set $A$ of $k-4$ colours such that 
$|F[p,q) \setminus A| = 1$ and  $|F[q,p) \setminus A| = 1$.
\item
$|F[1,l]| = k-1$ and there exist indices $1 \le p < q < r \le l$ and a set $A$ of $k-4$ colours such that 
$|F[p,q) \setminus A| = |F[q,r) \setminus A| = |F[r,p) \setminus A| = 1$.
\end{enumerate}

\begin{lemma}
\label{5color}
A $5$-colouring $f$ of a cycle of length $l$ is not bad if and only if it is good.
\end{lemma}

\noindent {\bf Proof:} We show that if $f$ is good it cannot satisfy any of the conditions in 
Lemma~\ref{unsafe} and is therefore not bad. This argument holds for all $k \ge 4$, but we state it
only for $k = 5$. Clearly $|F[1,l]| < 5$, hence $f$ does not satisfy condition 1 in 
Lemma~\ref{unsafe}. Similarly if $|F[1,l]|= 2$, $f$ does not satisfy any of the conditions in 
Lemma~\ref{unsafe}. Suppose $|F[1,l]|=3$ and $f$ satisfies condition 2 in the definition of a good
$5$-colouring with indices $p_1,q_1$. Then $f$ cannot satisfy condition 2 in the statement of 
Lemma~\ref{unsafe}. Suppose for contradiction $f$ satisfies condition 3 in Lemma~\ref{unsafe} with indices 
$p_2,q_2,r_2$. Since $[p_1,q_1), [q_1,p_1)$ is a partition of $[1,l]$, two of the indices $p_2,q_2,r_2$ 
are contained in one of the sets $[p_1,q_1)$ or $[q_1,p_1)$. Without loss of generality, assume $p_2,q_2 
\in [p_1,q_1)$. However, this implies $3 \le |F[p_2,q_2]| \le |F[p_1,q_1)| \le  2$, a contradiction. So $f$ 
cannot satisfy any of the conditions in Lemma~\ref{unsafe} and is not bad. 

Suppose $|F[1,l]| = 4$ and $f$ satisfies condition 3 in the definition of a good $5$-colouring with indices 
$p_1,q_1,r_1$ and let $A = \{c_1\}$. Suppose for contradiction $f$ satisfies condition 2 in Lemma~\ref{unsafe}
with indices $p_2,q_2$. Again, we may assume, without loss of generality, $p_1,q_1 \in [p_2,q_2)$. If $p_1 
\neq p_2$, then $[q_2,p_2] \subseteq [q_1,r_1) \cup [r_1,p_1)$, which contradicts the fact that 
$|F[q_2,p_2]| \ge 4$ but $|F[q_1,r_1) \cup F[r_1,p_1)| \le 3$. If $p_1 = p_2$ and $r_1 \in [p_2,q_2]$ then 
$[q_2,p_2] \subseteq [r_1,p_1]$, which is again a contradiction. If $r_1 \not\in [p_2,q_2]$ then $[p_2,q_2] 
\subseteq [p_1,q_1) \cup [q_1,r_1)$, which again contradicts the fact that $|F[p_2,q_2]| \ge 4$ and 
$|F[p_1,q_1) \cup F[q_1,r_1)| \le 3$. Finally, suppose $f$ satisfies condition 3 in Lemma~\ref{unsafe} 
with indices $p_2,q_2,r_2$. Let $c_2$ be a colour other than $c_1$ that is contained in $F[p_2,q_2] 
\cap F[q_2,r_2] \cap F[r_2,p_2]$. Without loss of generality, we may assume any vertex $v_m$ such that
$f(v_m) = c_2$ satisfies $m \in [p_1,q_1)$. However, this implies one of the sets $[p_2,q_2], [q_2,r_2], 
[r_2,p_2]$ must be contained in $[p_1,q_1)$ contradicting the fact that each of $F[p_2,q_2], F[q_2,r_2], 
F[r_2,p_2]$ has at least three elements, but $|F[p_1,q_1)| \le 2$. This proves that if $f$ is good, it is 
not bad.

We next show that if $f$ is not good then it is bad. This argument holds only for $k = 5$. We may assume 
$3 \le |F[1,l]| \le 4$, otherwise $f$ satisfies condition 1 in Lemma~\ref{unsafe} and is bad.

Suppose $|F[1,l]| = 3$. Let $p,q$ be indices such that $|[p,q]|$ is maximum and $|F[p,q)| = 2$.
Then $f(v_q) \neq f(v_m)$ and $f(v_{p-1}) \neq f(v_m)$ for any $m \in [p,q)$. Since $|F[1,l]| = 3$,
we must have $f(v_q) = f(v_{p-1})$ and $|F[q,p)| > 2$, otherwise $f$ satisfies condition 2 
in the definition of a good colouring with indices $p,q$. This implies $q \neq p-1$ and $[q,p)$ has at 
least four elements. Since $|F[p,q)| = 2$, we must have $q \ge p+2$. If $|F(p,q]| \ge 3$ then 
$\{p-1,p+1,q\}$ are three indices that satisfy condition 3 in Lemma~\ref{unsafe} and $f$ is bad.  If 
$|F(p,q]| = 2$, we must have $q = p+2$. The choice of $[p,q]$ then implies that $f(v_{q+1}) = f(v_p)$
and $f(v_{q+2}) = f(v_{p+1})$. This implies $f$ satisfies condition 3 in Lemma~\ref{unsafe} with
indices $\{p,q,q+2\}$ and is bad.

Suppose $|F[1,l]| = 4$. Let $p,q$ be indices such that $|F[p,q]| \ge 4$ and $|[p,q]|$ is minimum 
possible. We have $f(v_p) \neq f(v_m)$ for any $m \in (p,q]$, $f(v_q) \neq f(v_m)$ for any $m \in
[p,q)$, and also $f(v_p) \neq f(v_q)$.  This implies $|F(p,q)| = 2$. If $|F[q,p]| \ge 4$ then $p,q$ 
are indices satisfying condition 2 in Lemma~\ref{unsafe} and $f$ is bad. So $|F[q,p]| \le 3$,  
which implies $|F(p,q) \cap F[q,p]| \le 1$  and let $c_1$ be the colour in this set, if it is not empty,
otherwise let $c_1$ be any colour in $F(p,q)$. Suppose there exists an index $r \in [q,p)$ such that 
$f(v_p) \not\in F[q,r]$ and $f(v_q) \not\in F(r,p]$. Then $F(p,q), F[q,r], F(r,p]$ are sets of 
size at most 2, and each contains one element other than $c_1$. This implies $f$ satisfies condition 
3 in the definition of a good colouring, a contradiction. Therefore there is no such index $r$ and we have 
$p \neq q+1$. This implies there exists an index $r$ in $(q,p)$ such that $f(v_r) = f(v_p)$ and an index 
$r'$ in $(r,p)$ such that $f(v_{r'}) = f(v_q)$. If $q \ge p+4$ then we have $|F[p+2,r]| \ge 4$ and 
$|F[r,p+2]| \ge 4$, implying $f$ satisfies condition 2 in Lemma~\ref{unsafe} and is bad. 
Finally, suppose $q = p+3$. If either the colour $f(v_{p+1}) \in F[q,r']$ or $f(v_{p+2}) \in F[r,p]$, 
the same argument holds with indices $\{p+2,r'\}$ or $\{p+1,r\}$, respectively. Without loss of generality, 
assume the colour $c_1$ is $f(v_{p+1})$  so $v_{p+2}$ is the only vertex with the colour $f(v_{p+2})$. 
Choose $r'$ so that $|[r',p]|$ is as small as possible. Then $f(v_q) \not\in F(r',p]$ and $|F[p+2,q)| = 1$, 
$|F[q,r']| \le 2$ and $|F(r',p+2)| \le 2$. Also $|F[q,r'] \setminus \{f(v_p)\}| = 1$ and $|F(r',p+2)\setminus 
\{f(v_p)\}| = 1$, which implies the indices $p+2,q,r'+1$ with $A = \{f(v_p)\}$, satisfy condition 3 
in the definition of a good colouring, a contradiction. 
\hfill $\Box$

The $6$-colouring of the 8 cycle defined by $f(v_i) = f(v_{i+4}) = c_i$ for $1 \le i \le 4$ is an
example of a colouring that is neither good nor bad.

We can now state the main result.
\begin{theorem}

\label{main}
Any good $k$-colouring of a cycle is safe for all $k \ge 5$.
\end{theorem}

In order to prove Theorem~\ref{main} by induction, we need to consider sets of $k$-colourings of
a cycle, rather than a single $k$-colouring. We also need to consider near-triangulations $G$ with 
a given external boundary $C$ that may have chords. Let $C = v_1,v_2,\ldots,v_l$ be a cycle  and 
$L(v_i) \subseteq \{c_1,c_2,\ldots,c_k\}$  a non-empty list of at most $k$ colours assigned to the 
vertex $v_i$, for $1 \le i \le l$. A near-triangulation $G$ with external boundary $C$ is said to be 
\emph{consistent} with the list assignment if for any two vertices $u,v \in V(C)$ such that $uv$ is 
an edge or chord of $C$, both $u$ and $v$ are not assigned the same list of size 1. We will only
consider list assignments such that there are no two adjacent vertices $u,v$ in $C$ such that 
$|L(u)| = 1$ and $L(u) = L(v)$. We show that for certain kinds of list assignments to the vertices of $C$,
for any near-triangulation $G$ with external boundary $C$ that is consistent with it, there exists
an $L$-colouring of $C$ that can be extended to a $k$-colouring of $G$. 

Let $S_i = \{c_i\} \cup \{c_j | 5 \le j \le k\}$ for $i \in \{1,2,3\}$. A vertex $v$ is said to be of 
type $T_i$ for $i \in \{1,2,3\}$ if either $L(v) = \{c_1,c_2,c_3,c_4\} \setminus \{c_i\}$ or 
$|L(v)| = 1$ and $L(v) \subset S_i$. The list assignments that we consider will be such that every vertex 
will be of one of the three types.

\begin{lemma}
\label{one}
Let $C = v_1,v_2,\ldots,v_l$ be a cycle of length $l$ and $L$ an assignment of list of colours such that
every vertex is of type $T_1$. Let $(s,t) \in L(v_1) \times L(v_l)$ be any ordered pair of distinct
colours. Then for any near-triangulation $G$ with boundary $C$ that is consistent with the list
assignment, there exists an $L$-colouring $f$ of $C$ that extends to a $k$-colouring of $G$ with
$f(v_1) = s$ and $f(v_l) = t$.
\end{lemma}

\noindent {\bf Proof:} Suppose there is a counterexample and let $G$ be one with minimum number of edges.

\noindent{\bf Case 1.} Suppose $G$ has a chord $v_iv_j$ for some $i < j$. Then both $(i,j)$ and $(j,i)$ are 
non-empty sets. Let $C_1$ be the cycle $v_1,v_2,\ldots,v_i,v_j,v_{j+1},\ldots,v_l$ and $G_1$ the 
near-triangulation induced by the vertices on or in the interior of $C_1$. Then $G_1$ is consistent with the 
list assignment $L$ restricted to the vertices of $C_1$, and by the minimality of $G$, there exists an 
$L$-colouring $f_1$ of $C_1$ that can be extended to a $k$-colouring of $G_1$ with $f_1(v_1) = s$ and 
$f_1(v_l) = t$.  Let $C_2$ be the cycle $v_i,v_{i+1},\ldots, v_j$  and let $G_2$ be the near-triangulation 
induced by the vertices on or in the interior of $C_2$. Again, the list assignment to the vertices of $C_2$ 
satisfies the induction hypothesis and $G_2$ is consistent with it. The minimality of $G$ implies there exists 
an $L$-colouring $f_2$ of $C_2$ that can be extended to a $k$-colouring of $G_2$ with $f_2(v_i) = f_1(v_i)$ 
and $f_2(v_j) = f_1(v_j)$.  Defining $f(v) = f_1(v)$ for all vertices $v \in V(G_1)$ and $f(v) = f_2(v)$ for all 
vertices $v \in V(G_2)$, gives a $k$-colouring of $G$ that extends an $L$-colouring of $C$, a contradiction.

\noindent {\bf Case 2.} A similar argument holds if $G$ has a separating triangle $uvw$, that is a triangle 
whose interior as well as exterior contain vertices of $G$. Let $G_1$ be the near-triangulation obtained by 
deleting the vertices in the interior of the separating triangle $uvw$. Then $G_1$ has the same boundary $C$ 
as $G$ and is consistent with the list assignment $L$. The minimality of $G$ implies there exists an
$L$-colouring $f_1$ of $C$ that extends to a $k$-colouring of $G_1$, with $f_1(v_1) = s$ and $f_1(v_l) = t$. 
Rename the colours so that $f_1(u) = c_2$, $f_1(v) = c_3$ and $f_1(w) = c_4$. Let $G_2$ be the triangulation
induced by the vertices on and in the interior of the triangle $uvw$ and let $L(u) = L(v) = L(w) =
\{c_2,c_3,c_4\}$. Then the triangle $uvw$ satisfies the induction hypothesis and by induction, there exists
an $L$-colouring $f_2$ of the triangle that extends to a $k$-colouring of $G_2$, with $f_2(u) = c_2$ and
$f_2(v) = c_3$. The list assigned to $w$ ensures $f_2(w) = c_4$. Rename the colours in the colouring $f_2$ so 
that $f_2(x) = f_1(x)$ for all $x \in \{u,v,w\}$.  Again, setting $f(v) = f_1(v)$ for all vertices 
$v \in V(G_1)$ and $f(v)= f_2(v)$ for all vertices $v \in V(G_2)$, gives a $k$-colouring of $G$ that extends an
$L$-colouring of $C$, a contradiction.

\noindent {\bf Case 3.} We may now assume that $G$ is chordless and has no separating triangle.
If $G$ is a triangle, let $f(v_1) = s$ and $f(v_3) = t$ with $s \neq t$ by assumption. If $|L(v_1)| = 3$ then 
$s \in \{c_2,c_3,c_4\}$ and the same holds for $v_3$. This implies that if $|L(v_2)| = 1$, then 
$s,t \not\in L(v_2)$, and otherwise $|L(v_2)| = 3$ and $L(v_2)$ contains a colour other than $s,t$. 
Assigning this colour to $v_2$ gives an $L$-colouring of $C$.

Suppose $G$ is not a triangle, which implies that for every edge $v_iv_{i+1}$ in $C$, there
exists an internal vertex $v$ in $G$ such that $vv_iv_{i+1}$ is an internal face in $G$.
Suppose there exists such an edge $v_iv_{i+1}$ for some $1 \le i < l$ such that $L(v_i) \cap L(v_{i+1}) 
= \emptyset$. Let $G_1$ be obtained from $G$ by deleting the edge $v_iv_{i+1}$. Then $G_1$ is a 
near-triangulation bounded by the cycle $C_1 = v_1,\ldots,v_i,v,v_{i+1},\ldots,v_l$. Setting
$L(v) = \{c_2,c_3,c_4\}$ gives a list assignment to the vertices of $C_1$ that satisfies the induction 
hypothesis, and $G_1$ is consistent with it. The minimality of $G$ implies there
exists an $L$-colouring  $f$ of $C_1$ that extends to a $k$-colouring of $G_1$ with $f(v_1) = s$
and $f(v_l)=t$. Since $L(v_i) \cap L(v_{i+1}) = \emptyset$ by assumption, this gives the required
$L$-colouring of $C$ that extends to a $k$-colouring of $G$.  If there is no such edge $v_iv_{i+1}$ in $C$,
then for all $m \in (1,l)$, $L(v_m) = \{c_2,c_3,c_4\}$, for otherwise $L(v_i) = \{a\}$ for some 
$a \in S_1$ and  $a \not\in L(v_{i+1})$, which implies $L(v_i) \cap L(v_{i+1}) = \emptyset$. This also 
implies $L(v_1) = L(v_l) = \{c_2,c_3,c_4\}$.  

Let $v_1 = u_1, u_2,\ldots, u_r = v_3$ be the vertices adjacent to $v_2$ such that $v_2u_iu_{i+1}$ 
is an internal face of $G$, for some $r \ge 3$ and all $1 \le i < r$. Let $G_1$ be the graph obtained 
from $G$ by deleting the vertex $v_2$. Then $G_1$ is a near-triangulation bounded by the cycle $C_1 = 
v_1,u_2,\ldots,u_{r-1},v_3,\ldots,v_l$. Let $L(u_i) = \{c_1\}$ if $i$ is odd and $L(u_i) = \{c_5\}$ if 
$i$ is even, for all $1 < i < r$. The list assignment $L$ to the vertices of $C_1$ satisfies the induction 
hypothesis and $G_1$ is consistent with it. The minimality of $G$ implies there exists an $L$-colouring $f$ of 
$C_1$ that extends to a $k$-colouring of $G_1$ with $f(v_1) = s$ and $f(v_l) = t$. Since $L(v_2) = 
\{c_2,c_3,c_4\}$, setting $f(v_2)$ to be a colour in $\{c_2,c_3,c_4\} \setminus (\{f(v_1)\} \cup \{f(v_3)\})$, 
gives  the required $L$-colouring of $C$ that extends to a $k$-colouring of $G$. This completes all cases 
and the proof of the lemma.
\hfill $\Box$

We next consider the case when the vertices can be of two types. In this case, we need to put
more restrictions on the $L$-colouring of the cycle. If $uv$ is an edge in the cycle $C$,
we place restrictions on the ordered pair of colours $(f(u),f(v))$ assigned to the
vertices $u$ and $v$ in an $L$-colouring $f$. In the usual list colouring, this can be any
pair in $L(u) \times L(v)$ such that the two colours are distinct. Here we allow only specified subsets of 
such pairs and assign a list of allowed pairs from $L(u) \times L(v)$ to the edge $uv$. We call these sets 
of ordered pairs the labels of the edges. Let $a$ denote any colour in $S_1$, $b$ denote any colour in $S_2$,  
and $x$ denote any colour in $\{c_3,c_4\}$ with $\{x,y\} = \{c_3,c_4\}$. The allowed label sets are 
the following and we call the collection of these $\mathcal{L}_{12}$.

\begin{tabular}{lll}
(i) $\{(a,b)\}, a \neq b$ & (ii) $\{(a,x)\}$ & (iii)  $\{(x,b)\}$ \\
(iv) $\{(c_2,x),(y,x)\}$ & (v) $\{(x,c_1),(x,y)\}$ & (vi) $\{(c_2,x), (x,c_1)\}$ 
\end{tabular}

We consider list assignments to a cycle $C$ in which every vertex $v$ is assigned a list of colours and
exactly two edges $u_1v_1$ and $u_2v_2$ are assigned labels $L(u_1v_1), L(u_2v_2) \in \mathcal{L}_{12}$ 
respectively, where $L(u_1v_1) \subseteq L(u_1) \times L(v_1)$ and $L(u_2v_2) \subseteq L(u_2) \times L(v_2)$. 
Let $G$ be a near-triangulation with boundary $C$ that is consistent with the lists assigned to the vertices 
in $C$. We say the list assignment is \emph{feasible} for $G$ if for every ordered pair $(s_1,t_1) \in 
L(u_1v_1)$ there exists an ordered pair $(s_2,t_2) \in L(u_2v_2)$ and an $L$-colouring $f$ of $C$ that can be 
extended to a $k$-colouring of $G$ with $f(u_i) = s_i$ and $f(v_i) = t_i$ for $i \in \{1,2\}$. We call
$(L(u_1v_1),L(u_2v_2))$ a feasible pair of labels for $G$.

\begin{lemma}
\label{two}
Let $C = v_1,v_2,\ldots,v_l$ be a cycle of length $l$ and $p$ an integer such that $1 \le p < l$.
Let $L$ be an assignment of  list of colours to each vertex in $C$ such that $v_m$ is of type $T_1$
for all $m \in [1,p]$, and of type $T_2$ for all $m \in (p,1)$. Let $L_2=L(v_pv_{p+1}) \subseteq L(v_p) \times 
L(v_{p+1})$ be any label in $\mathcal{L}_{12}$ assigned to the edge $v_pv_{p+1}$. For any near-triangulation $G$ 
with boundary $C$ that is consistent with $L$, there exists a label $L_1 \in \mathcal{L}_{12}$ and
$L_1  \subseteq L(v_1) \times L(v_l)$ that can be assigned to the edge $v_1v_l$ such that $(L_1, L_2)$ is a 
feasible pair of labels for $G$.
\end{lemma} 

\noindent {\bf Proof:} The proof is again by induction on the number of edges. Let $G$ be a counterexample
with the minimum number of edges. 

\noindent {\bf Case 1.} Suppose $G$ has a chord $v_iv_j$ with $1 \le i < j \le l$.  Let $C_1$ be 
the cycle $v_1,v_2,\ldots,v_i,v_j,\ldots,v_l$ and $G_1$ the near-triangulation induced by the vertices on and 
in the interior of $C_1$. Let $G_2$ be the near-triangulation induced by the vertices on and in the interior 
of the cycle $C_2 = v_i,v_{i+1},\ldots,v_j$. 

\noindent {\bf Case 1.1} Suppose $j \le p$. The list assignment $L$ restricted to the vertices in $C_1$, along 
with the label $L_2$ assigned to the edge $v_pv_{p+1}$, satisfies the induction hypothesis, and the minimality 
of $G$ implies the existence of a feasible pair of labels $(L_1,L_2)$ for $G_1$. We claim that this is also 
feasible for $G$. Consider any $L$-colouring $f_1$ of $C_1$ that can be extended to a $k$-colouring of $G_1$. 
The list assignment $L$ restricted to the vertices of $C_2$ satisfies the conditions of Lemma~\ref{one} and 
applying it gives an $L$-colouring $f_2$ of $C_2$ that can be extended to a $k$-colouring of $G_2$, with 
$f_2(v) = f_1(v)$ for $v \in \{v_i,v_j\}$. Thus any $L$-colouring of $C_1$ that can be extended to a 
$k$-colouring of $G_1$ can also be extended to a $k$-colouring of $G$. Therefore the labels $(L_1,L_2)$ form a 
feasible pair for $G$.  A similar argument holds if $i > p$. In this case, for applying Lemma~\ref{one} to the 
near-triangulation $G_2$, first swap the colours $c_1$ and $c_2$ in the lists assigned to the vertices of $C_2$. 
This converts all vertices $v_m$ with $m \in (i,j)$ to type $T_1$ and Lemma~\ref{one} can be applied. After 
getting the $L$-colouring of $C_2$ that can be extended to a $k$-colouring of $G_2$, again swap the colours 
$c_1,c_2$ in the $k$-colouring of $G_2$ to get the required $L$-colouring of $C_2$.

\noindent {\bf Case 1.2} Suppose $i \le p < j$. In this case, the list assignment $L$ restricted to the vertices 
of $C_2$ satisfies the induction hypothesis, considering $v_iv_j$ to be the edge $v_1v_l$.  The minimality of 
$G$ implies there exists a label $L_1'$ that can be assigned to the edge $v_iv_j$ such that $(L_1',L_2)$ is a 
feasible pair of labels for $G_2$. The list assignment $L$ restricted to the vertices of $C_1$, with the 
label $L_1'$ assigned to the edge $v_iv_j$, also satisfies the induction hypothesis, and there exists a 
label $L_1$ that can be assigned to the edge $v_1v_l$ such that $(L_1,L_1')$ is a feasible pair for $G_1$. 
This implies $(L_1,L_2)$ is a feasible pair of labels for $G$.

\noindent {\bf Case 2.} We may assume $G$ is chordless. If $G$ has a separating triangle, the same 
argument as in Lemma~\ref{one} can be used. Delete the vertices in the interior of a separating triangle to
get a near-triangulation $G_1$ with the same boundary as $G$. By induction, there exists a feasible
pair of labels $(L_1,L_2)$ for $G_1$. It can then be argued as in Lemma~\ref{one} that $(L_1,L_2)$ is 
also feasible for $G$. We may therefore assume $G$ has no chord and no separating triangle. 

\noindent {\bf Case 2.1} Suppose $l=3$ and without loss of generality assume $p=2$, otherwise reverse the cycle 
and swap the colours $c_1$ and $c_2$ to apply the argument. In this case, it is sufficient to show that for 
every label $L_2 \in \mathcal{L}_{12}$ assigned to the edge $v_2v_3$, there exists a label $L_1$ that can be 
assigned to $v_1v_3$ such that for every pair $(s,t) \in L_1$, there exists a colour $r \not\in \{s,t\}$ such 
that $(r,t) \in L_2$. Lemma~\ref{one} implies that the $L$-colouring $f$ of $C$ defined by $f(v_1) = s$, 
$f(v_2) = r$, $f(v_3)=t$, can be extended to a $k$-colouring of $G$, hence $(L_1,L_2)$ is a feasible pair for 
$G$.

If $|L(v_2)| = |L(v_3)| = 1$ then $L_2 = L(v_2) \times L(v_3)$. If $|L(v_1)| = 1$, then since $G$ is consistent 
with $L$, $L(v_1) \neq L(v_2)$ and $L(v_1) \neq L(v_3)$. Setting $L_1 = L(v_1) \times L(v_3)$ gives a feasible 
pair of labels $(L_1,L_2)$ for $G$. If $|L(v_1)| = 3$, then since $v_1$ is of type $T_1$, $L(v_1) = 
\{c_2,c_3,c_4\}$. Setting $L_1 = \{c_3\} \times L(v_3)$ gives the required label for the edge $v_1v_3$. 

Suppose $|L(v_2)|=3$ and $|L(v_3)| = 1$. We can assume $L_2 = \{(x,b)\}$ for some $x \in \{c_3,c_4\}$ 
and $b \in S_2$. If $|L(v_1)| = 1$, let $L_1$ be the set $L(v_1) \times \{b\}$. If $|L(v_1)| = 3$, let $L_1$ 
be the set $\{(y,b)\}$, where $y \in \{c_3,c_4\} \setminus \{x\}$. In both cases, we get a feasible pair of 
labels $(L_1,L_2)$ for $G$.

Suppose $|L(v_2)| = 1$ and $|L(v_3)| = 3$. We can assume that $L_2 = \{(a,x)\}$ for some $x \in 
\{c_3,c_4\}$ and $a \in S_1$. If $|L(v_1)| = 1$ let $L_1 = L(v_1) \times \{x\}$, otherwise $L(v_1) = 
\{c_2,c_3,c_4\}$ and let $L_1 = \{(c_2,x),(y,x)\}$. In both cases $(L_1,L_2)$ is a feasible pair of labels 
for $G$.

Suppose both $|L(v_2)|$ and $|L(v_3)|$ are 3, which implies $L(v_2) = \{c_2,c_3,c_4\}$ and
$L(v_3) = \{c_1,c_3,c_4\}$. Then $L_2$ can be any one of $\{(c_2,x),(y,x)\}$, $\{(x,c_1),(x,y)\}$
and $\{(c_2,x),(x,c_1)\}$. If $|L(v_1)| = 1$, then $L(v_1) = \{a\}$ for some $a \in S_1$. If $(c_2,x) 
\in L_2$, let $L_1 = \{(a,x)\}$, otherwise $(x,y) \in L_2$ and let $L_1 = \{(a,y)\}$. In both cases 
$(L_1,L_2)$ is a feasible pair for $G$. If $|L(v_1)| = 3$, then $L(v_1) = \{c_2,c_3,c_4\}$. If 
$L_2 = \{(c_2,x),(y,x)\}$ let $L_1 = L_2$. If $L_2 = \{(x,c_1),(x,y)\}$ let $L_1 = \{(y,c_1),(c_2,y)\}$ and 
if $L_2 = \{(c_2,x),(x,c_1)\}$ let $L_1 = \{(y,x),(y,c_1)\}$. In each case, it is easy to verify that
$(L_1,L_2)$ is a feasible pair of labels for $G$. 

\noindent {\bf Case 2.2} Suppose now that $l > 3$. Since $G$ is chordless, for every edge $v_iv_{i+1}$ in $C$, 
there exists an internal vertex $v$ such that $vv_iv_{i+1}$ is an internal face in $G$. Suppose there exists 
such an edge $v_iv_{i+1}$ for some $1 \le i < l$ and $i \neq p$ such that $L(v_i) \cap L(v_{i+1}) = \emptyset$. 
Let $G_1$ be the near-triangulation obtained by deleting the edge $v_iv_{i+1}$, with the cycle $C_1 = v_1, v_2, 
\ldots, v_i, v, v_{i+1},\ldots,v_l$ as the boundary. If $i < p$, let $L(v) = \{c_2,c_3,c_4\}$ and if $i > p$ 
then let $L(v) = \{c_1,c_3,c_4\}$. The list assignment $L$ to the vertices of $C_1$, along with the label $L_2$ 
to the edge $v_pv_{p+1}$, satisfies the induction hypothesis and $G_1$ is consistent with it. By induction, 
there exists a label $L_1$ for the edge $v_1v_l$ such that $(L_1,L_2)$ is a feasible pair of labels for $G_1$. 
Since every $L$-colouring of $C_1$ is also an $L$-colouring of $C$, $(L_1,L_2)$ is a feasible pair 
of labels for $G$ also. We may assume there is no such edge in $C$. 

\noindent {\bf Case 2.2.1} Suppose there exists a vertex $v_i$ with $|L(v_i)| = 1$. Suppose $i = 1$,
which implies $L(v_1) = \{a\}$ for some $a \in S_1$. If $p > 1$, then $L(v_1) \cap L(v_2) = \emptyset$,
a contradiction. Therefore $p = 1$ and the label $L_2$ attached to the edge $v_1v_2$ is either $\{(a,b)\}$ for
some $b \in S_2$ or $\{(a,x)\}$ for $x \in \{c_3,c_4\}$. Let $v$ be the internal vertex in $G$ such that
$vv_1v_2$ is an internal face in $G$. Let $G_1$ be the near-triangulation obtained by deleting the edge
$v_1v_2$  with the cycle $C_1 = v_1,v,v_2,\ldots,v_l$ as the boundary. Let $L(v) = \{c_2,c_3,c_4\}$
and assign the label $L_2'$ to the edge $vv_2$, where $L_2' = \{(c_3,b)\}$ if $L_2 = \{(a,b)\}$ and
$L_2' = \{(c_2,x),(y,x)\}$ if $L_2 = \{(a,x)\}$. The list assignment $L$ to the vertices in $C_1$, along
with the label $L_2'$ assigned to the edge $vv_2$, satisfies the induction hypothesis, and by induction,
there exists a feasible pair $(L_1,L_2')$ for $G_1$. Since any $L$-colouring of $C_1$ is also an
$L$-colouring of $C$, $(L_1,L_2)$ is a feasible pair for $G$. A symmetrical argument holds if $|L(v_2)|=1$,
by reversing the cycle and swapping the colours $c_1,c_2$.

If $1 < i \le p$ then $L(v_{i-1}) \cap L(v_i) = \emptyset$ and if $l > i > p$ then $L(v_i) \cap L(v_{i+1}) 
= \emptyset$, a contradiction.
 
\noindent {\bf Case 2.2.2} Suppose $|L(v_i)| = 3$ for all $i \in [1,l]$ and suppose $p > 2$. Let $v_1 = 
u_1,u_2,\ldots,u_r = v_3$ be the vertices adjacent to $v_2$ such that $v_2u_iu_{i+1}$ is an internal face in $G$ 
for all $1 \le i < r$. Let $G_1$ be obtained from $G$ by deleting the vertex $v_2$. Then $G_1$ is a 
near-triangulation bounded by the cycle $C_1 = v_1, u_2, \ldots,u_{r-1}, v_3, \ldots, v_l$. Let $L(u_i) = 
\{c_1\}$ if $i$ is odd and $L(u_i) = \{c_5\}$ if $i$ is even, for all $2 \le i < r$. This gives a list 
assignment to the vertices of $C_1$, which along with the label $L_2$ for the edge $v_pv_{p+1}$,  satisfies the 
induction hypothesis and $G_1$ is consistent with it. By induction, there exists a label $L_1$ for the edge 
$v_1v_l$ such that $(L_1,L_2)$ is a feasible pair of labels for $G_1$. Since any $L$-colouring of $C_1$ can be 
extended to an $L$-colouring of $C$ by assigning $v_2$ a colour in $\{c_2,c_3,c_4\}$ that is not assigned to 
$v_1$ or $v_3$, $(L_1,L_2)$ is also a feasible pair of labels for $G$. A symmetrical argument can be used if 
$p < l-2$, by relabeling the vertex $v_i$ as $v_{l+1-i}$ and swapping the colours $c_1$ and $c_2$.

\noindent {\bf Case 2.2.3} Suppose $|L(v_i)| = 3$ for all $i \in [1,l]$, $p \le 2$ and $p \ge l-2$, which 
implies $l = 4$, $p = 2$. This implies $L_2$ is one of the labels $\{(c_2,x),(y,x)\}$, $\{(x,c_1), (x,y)\}$ or 
$\{(c_2,x), (x,c_1)\}$. 

Suppose $(x,c_1) \in L_2$. Again consider the near-triangulation $G_1$ obtained by deleting the vertex $v_2$
and bounded by the cycle $C_1 = v_1, u_2,\ldots,v_3,v_4$. Let $L(u_r) = \{c_1\}$ and $L(u_i) = \{c_1,c_5\}
\setminus L(u_{i+1})$ for $2 \le i < r$.  Note that the list for $v_3 = u_r$ is modified by removing the 
elements $c_3,c_4$ from it. Assign the label $L_2' = \{(c_1,x)\}$ to the edge $v_3v_4$. Then the list 
assignment to the vertices of $C_1$ satisfies the induction hypothesis, and $G_1$ is consistent with it. By 
induction, there exists a label $L_1$ for the edge $v_1v_4$ such that $(L_1,L_2')$ is a feasible pair of labels 
for $G_1$. Since in any $L$-colouring $f$ of $C_1$ we must have  $f(v_4) = x$,  $L_1 = \{(c_2,x),(y,x)\}$. 
Therefore any $L$-colouring of $C_1$ can be extended to an $L$-colouring of $C$ by assigning colour $x$ to 
$v_2$. This implies $(L_1,L_2)$ is a feasible pair of labels for $G$. If $(x,c_1) \not\in L_2$, we can assume
$(c_2,x) \in L_2$, and a symmetrical argument can be used after reversing the cycle and swapping the
colours $c_1, c_2$.

This completes all cases and the proof of Lemma~\ref{two}.
\hfill $\Box$

We now consider list assignments with vertices of 3 different types. We define two other sets of possible
labels for edges in the cycle. Let $a \in S_1$, $b \in S_2$ and $c \in S_3$ be any colours. The label set 
$\mathcal{L}_{13}$ contains the following sets of ordered pairs.

\begin{tabular}{lll}
(i) $\{(a,c)\}$ & (ii) $\{(a,c_2), (a,c_4)\}$ & (iii) $\{(c_2,c), (c_4,c)\}$ 
\end{tabular}

\begin{tabular}{l}
(iv) $\{(c_3,c_2), (c_3,c_4), (c_2,c_4), (c_4,c_2)\}$ \\
(v) $\{(c_2,c_1),(c_4,c_1), (c_2,c_4), (c_4,c_2) \}$ \\
(vi) $\{(c_3,c_2),(c_3,c_4),(c_2,c_1),(c_4,c_1)\}$
\end{tabular}
 
A set in $\mathcal{L}_{13}$ is obtained from a set in $\mathcal{L}_{12}$ by first adding all pairs obtained by 
swapping colours $c_3$ and $c_4$ and then swapping the colours $c_2$ and $c_3$. Thus if the set in 
$\mathcal{L}_{12}$ is  $\{(c_2,c_3), (c_3,c_1)\}$, adding all pairs obtained by swapping $c_3$ and $c_4$ gives 
the set $\{(c_2,c_3), (c_3,c_1), (c_2,c_4), (c_4,c_1)\}$ and then swapping the colours  $c_2$ and $c_3$ gives 
the label $\{(c_3,c_2), (c_2,c_1), (c_3,c_4), (c_4,c_1)\} \in \mathcal{L}_{13}$. The labels in 
$\mathcal{L}_{32}$ are obtained in a similar way from those in $\mathcal{L}_{12}$, by first adding all pairs 
obtained by swapping colours $c_3$ and $c_4$ and then swapping the colours $c_1$ and $c_3$. For the example 
given, the label obtained is $\{(c_2,c_1), (c_1,c_3), (c_2,c_4), (c_4,c_3)\}$. The label set $\mathcal{L}_{32}$ 
contains the following sets of ordered pairs.

\begin{tabular}{lll} 
(i) $\{(c,b)\}$ & (ii) $\{(c,c_1), (c,c_4)\}$ & (iii) $\{(c_1,b), (c_4,b)\}$
\end{tabular}

\begin{tabular}{l} 
(iv) $\{(c_2,c_1), (c_2,c_4), (c_1,c_4),(c_4,c_1)\}$ \\
(v)  $\{(c_1,c_3), (c_4,c_3), (c_1,c_4), (c_4,c_1)\}$ \\
(vi) $\{ (c_2,c_1), (c_2,c_4), (c_1,c_3), (c_4,c_3) \}$. 
\end{tabular}

The sets in $\mathcal{L}_{13}$ and $\mathcal{L}_{32}$ correspond to the sets in
$\mathcal{L}_{12}$ from which they are obtained.

Note that in Lemma~\ref{two}, the list assignments and labels are symmetric in the colours $c_3$ and
$c_4$. Therefore if $(L_1,L_2)$ is a feasible pair for a near-triangulation $G$, then so is
$(L_1',L_2')$, where $L_1'$ and $L_2'$ are obtained from $L_1$ and $L_2$ by swapping the colours
$c_3$ and $c_4$. 

\begin{lemma}
\label{three}
Let $C= v_1,v_2,\ldots,v_l$ be a cycle of length $l$ and $p,q$ positive integers such that $1 \le p < q < l$.
Let $L$ be an assignment of list of colours to the vertices of $C$ such that the vertex $v_m$ is of type
$T_1$ for all $m \in [1,p]$, of type $T_3$ for $m \in (p,q]$ and of type $T_2$ for $m \in (q,l]$. Suppose
$L_2 \in \mathcal{L}_{13}$ and $L_3 \in \mathcal{L}_{32}$ are labels assigned to the edges $v_pv_{p+1}$
and $v_qv_{q+1}$, respectively, such that $L_2 \subseteq L(v_p) \times L(v_{p+1})$ and $L_3 \subseteq
L(v_q) \times L(v_{q+1})$. Then for any near-triangulation $G$ with boundary $C$ that is consistent
with the list assignment $L$, there exists a label $L_1 \in \mathcal{L}_{12}$ and $L_1 \subseteq L(v_1)
\times L(v_l)$ that can be assigned to the edge $v_1v_l$, such that for any ordered pair $(s_1,t_1) \in L_1$, 
there exist ordered pairs $(s_2,t_2) \in L_2$, $(s_3,t_3) \in L_3$ and an $L$-colouring $f$ of $C$ that can be 
extended to a $k$-colouring of $G$ with $f(v_1) = s_1$, $f(v_l) = t_1$, $f(v_p) = s_2$, $f(v_{p+1}) = t_2$, 
$f(v_q) = s_3$ and $f(v_{q+1}) = t_3$.
\end{lemma}

\noindent {\bf Proof:} The proof is again by induction, and we consider a counterexample $G$ with
the minimum number of edges. We call the triple of labels $(L_1,L_2,L_3)$ satisfying the properties
in the lemma a feasible triple, and suppose $G$ does not have a feasible triple for some list
assignment $L$ to the vertices of $C$ and labels $L_2,L_3$.

\noindent {\bf Case 1.} Suppose $G$ has a chord $v_iv_j$ for some $1 \le i < j \le l$. Again, let $G_1$ be the 
near-triangulation induced by the vertices on or in the interior of the cycle $C_1 = v_1,\ldots,v_i,v_j,
\ldots,v_l$ and $G_2$ the near-triangulation induced by the vertices on or in the interior of the cycle 
$C_2 = v_i,v_{i+1}\ldots,v_j$.

\noindent {\bf Case 1.1} Suppose $v_i$ and $v_j$ are vertices of the same type. Then the list assignment $L$ 
restricted to the vertices of $C_1$, along with the labels $L_2$ and $L_3$ assigned to the edges $v_pv_{p+1}$ 
and $v_qv_{q+1}$, satisfies the induction hypothesis. By induction, there exists a feasible triple
$(L_1,L_2,L_3)$ for $G_1$, and we claim that it is also feasible for $G$. Let $f_1$ be any $L$-colouring
of $C_1$. The list assignment $L$ restricted to the vertices of $C_2$ satisfies the conditions of
Lemma~\ref{one}, perhaps after renaming colours so that vertices $v_m$ have type $T_1$ for $m \in [i,j]$.
Lemma~\ref{one} then implies there exists an $L$-colouring $f_2$ of $C_2$ that can be extended to a $k$-colouring
of $G_2$ with $f_2(v_i) = f_1(v_i)$ and $f_2(v_j) = f_1(v_j)$. Since this holds for any $L$-colouring of $C_1$, 
$(L_1,L_2,L_3)$ is a feasible triple for $G$.

\noindent {\bf Case 1.2} Suppose $v_i$ is of type $T_1$ and $v_j$ of type $T_2$, which implies $1 \le i \le p$ 
and $q < j \le l$. The list assignment $L$ restricted to the vertices of $C_2$, along with the labels $L_2,L_3$, 
satisfies the induction hypothesis, and by induction there exists a label $L_1' \in \mathcal{L}_{12}$ that can 
be assigned to the edge $v_iv_j$ such that $(L_1',L_2,L_3)$ is a feasible triple for $G_2$. The list assignment 
$L$ restricted to the vertices of $C_1$, along with the label $L_1'$ assigned to the edge $v_iv_j$, satisfies 
the conditions of Lemma~\ref{two}. Applying Lemma~\ref{two} to $G_1$, there exists a label $L_1$ that can be 
assigned to the edge $v_1v_l$ such that $(L_1,L_1')$ is a feasible pair of labels for $G_1$. This implies
$(L_1,L_2,L_3)$ is a feasible triple for $G$.

\noindent {\bf Case 1.3} Suppose $v_i$ is of type $T_1$ and $v_j$ of type $T_3$. Then the list assignment $L$ 
restricted to the cycle $C_1$ satisfies the induction hypothesis but the edge $v_pv_{p+1}$ is not in $C_1$. We 
choose an appropriate label $L_2' \in \mathcal{L}_{13}$ for the edge $v_iv_j$ and the label $L_3$ for the
edge $v_qv_{q+1}$ in $C_1$. Applying induction gives a feasible triple $(L_1,L_2',L_3)$ for $G_1$.
We choose $L_2'$ so that for any ordered pair $(s_2',t_2') \in L_2'$, there exists an ordered pair
$(s_2,t_2) \in L_2$ and an $L$-colouring $f_2$ of $C_2$ that can be extended to a $k$-colouring
of $G_2$ with $f_2(v_i) = s_2'$, $f_2(v_j) = t_2'$, $f_2(v_p) = s_2$ and $f_2(v_{p+1}) = t_2$.
Then $(L_1,L_2,L_3)$ is a feasible triple for $G$.

The label $L_2'$ for the edge $v_iv_j$ is found by applying Lemma~\ref{two} to the near-triangulation
$G_2$. We first swap the colours $c_2$ and $c_3$ in all the lists and the ordered pairs in the label
so that vertices $v_m$ of type $T_3$ for $m \in (p,j]$ are converted to type $T_2$. This swap does not affect 
the lists for the vertices $v_m$ of type $T_1$ for $m \in [i,p]$. After swapping colours $c_2$ and $c_3$,
the labels in $\mathcal{L}_{13}$ are closed under swapping colours $c_3$ and $c_4$. Retaining only one of the
pairs that can be obtained by such a swap converts the label $L_2$ into a label $L_2''$ that is in 
$\mathcal{L}_{12}$. In other words, $L_2''$ is the label in $\mathcal{L}_{12}$ corresponding to the label
$L_2 \in \mathcal{L}_{13}$.  Lemma~\ref{two} implies there exists a label $L_2''' \in \mathcal{L}_{12}$ that 
can be assigned to the edge $v_iv_j$ such that $(L_2''',L_2'')$ is a feasible pair for $G_2$. Let $L_2'$ be
the label in $\mathcal{L}_{13}$ corresponding to the label $L_2'''$. It follows from Lemma~\ref{two}, that 
$(L_2',L_2)$ is a feasible pair of labels for $G_2$. 

The argument in the case $v_i$ is of type $T_3$ and $v_j$ is of type $T_2$ is symmetric, and the
above argument can be applied after reversing the cycle and swapping the colours $c_1$ and
$c_2$ in the lists and the labels. 

\noindent {\bf Case 2.} We may now assume $G$ has no chords. If $G$ has a separating triangle, exactly the 
same argument as in Lemmas ~\ref{one} and ~\ref{two} can be used. So we may assume $G$ has no chords or 
separating triangles.

\noindent {\bf Case 2.1} Suppose $l=3$, which implies $p = 1$ and $q = 2$. In this case, it is sufficient to 
show that for any labels $L_2,L_3$ assigned to the edges $v_1v_2$ and $v_2v_3$, respectively, there exists a 
label $L_1$ that can be assigned to $v_1v_3$ such that for every pair $(s,t) \in L_1$, there is colour $r \not\in 
\{s,t\}$ that can be assigned to $v_2$ such that $(s,r) \in L_2$ and $(r,t) \in L_3$. Lemma~\ref{one} implies 
that any such $L$-colouring of $C$ can be extended to a $k$-colouring of $G$, hence $(L_1,L_2,L_3)$ is a 
feasible triple for any triangulation $G$.

Suppose the label $L_2$ is $\{(a,c)\}$ for some $a \in S_1$ and $c \in S_3$. Then $L_3$ is
either $\{(c,b)\}$ for some $b \in S_2$ or $\{(c,c_1),(c,c_4)\}$. In the first case, the $L$-colouring
is given by $f(v_1) = a$, $f(v_2) = c$ and $f(v_3) = b$, hence $L_1=\{(a,b)\}$. In the second case,
set $f(v_3) = c_4$ and hence $L_1 = \{(a,c_4)\}$ satisfies the requirements.

Suppose $L_2$ is the label $\{(a,c_2),(a,c_4)\}$ for some $a \in S_1$. If $L_3 = \{(c_1,b),(c_4,b)\}$,
for some $b \in S_2$, let $f$ be the $L$-colouring of $C$ with $f(v_1) = a$, $f(v_2) = c_4$ and 
$f(v_3) = b$. In this case, the label $L_1$ is $\{(a,b)\}$. Otherwise if $(c_1,c_3),(c_4,c_3) \in L_3$, 
set $f(v_3) = c_3$ to get the required label $L_1 = \{(a,c_3)\}$. The remaining possibility is
that $L_3 = \{(c_2,c_1),(c_2,c_4),(c_1,c_4),(c_4,c_1)\}$. In this case, choose $f(v_2) = c_2$ and
$f(v_3) = c_4$, to get the required label $L_1 = \{(a,c_4)\}$.

We may now assume $|L(v_1)| = 3$ and by symmetry, $|L(v_3)| = 3$. Suppose $|L(v_2)| = 1$, which
implies $L(v_2) = \{c\}$ for some $c \in S_3$. In this case $L_2 = \{(c_2,c),(c_4,c)\}$ and
$L_3 = \{(c,c_1),(c,c_4)\}$. Then setting $f(v_1) = c_2$, $f(v_2) = c$ and $f(v_3) = c_4$ gives
an $L$-colouring of $C$. Similarly, $f(v_1) = c_4$, $f(v_2) = c$ and $f(v_3) = c_1$ is an
$L$-colouring of $C$. This implies $L_1 = \{(c_2,c_4),(c_4,c_1)\}$ gives a feasible triple for
$G$.

Finally, suppose all three vertices have lists of size 3. Suppose $(c_2,c_1),(c_4,c_1) \in L_2$.
If $(c_1,c_3),(c_4,c_3) \in L_3$, let $f(v_2) = c_1$, $f(v_3) = c_3$ and $f(v_1)$ can be either
$c_2$ or $c_4$. These give $L$-colourings of $C$, hence setting $L_1 = \{(c_2,c_3),(c_4,c_3)\}$
gives a feasible triple for $G$. If $(c_1,c_3) \not\in L_3$ then $L_3 = \{(c_2,c_1),(c_2,c_4),
(c_1,c_4),(c_4,c_1)\}$. Applying a symmetrical argument by reversing the cycle and swapping the
colours $c_1,c_2$, we must have $L_2 = \{(c_2,c_1),(c_4,c_1), (c_2,c_4),(c_4,c_2)\}$. In this
case, setting $f(v_1) = c_2$, $f(v_2) = c_1$ and $f(v_3) = c_4$ gives an $L$-colouring of $C$,
as does setting $f(v_1) = c_4$, $f(v_2) = c_2$ and $f(v_3) = c_1$. This implies $L_1 =
\{(c_2,c_4),(c_4,c_1)\}$ gives a feasible triple for $G$. The final case to consider is if
$(c_2,c_1),(c_4,c_1) \not\in L_2$ and by a symmetrical argument $(c_2,c_1),(c_2,c_4) \not\in
L_3$. This implies $L_2 = \{(c_3,c_2),(c_3,c_4),(c_2,c_4),(c_4,c_2)\}$ and
$L_3 =  \{(c_1,c_3),(c_4,c_3),(c_1,c_4),(c_4,c_1)\}$.  Therefore $f(v_1) = c_3$, $f(v_2) =
c_4$ and $f(v_3) = c_1$ is an $L$-colouring of $C$ and so is $f(v_1) = c_2$, $f(v_2) = c_4$
and $f(v_3) = c_3$. This implies $L_1 = \{(c_2,c_3),(c_3,c_1)\}$ gives a feasible triple 
for $G$. This completes all cases when $l = 3$.

\noindent {\bf Case 2.2} We may now assume $l > 3$  and $G$ is a chordless near-triangulation with no 
separating triangles. Therefore for every edge $v_iv_{i+1}$ in $C$ there exists an internal vertex $v$ in
$G$ such that $vv_iv_{i+1}$ is an internal face in $G$. As in the proof of Lemma~\ref{two}, we
may assume that $L(v_i) \cap L(v_{i+1}) \neq \emptyset$ for any $i \not\in \{p,q\}$, otherwise we can get 
the required feasible triple of labels by deleting such an edge and applying induction.

\noindent {\bf Case 2.2.1} Suppose there exists a vertex $v_i$ with $|L(v_i)| = 1$. Suppose $i = 1$,
which implies $L(v_1) = \{a\}$ for some $a \in S_1$. Then we must have $p = 1$ otherwise $v_2$ is
also a vertex of type $T_1$ and $L(v_1) \cap L(v_2) = \emptyset$. This implies the label $L_2$ assigned
to the edge $v_1v_2$ is either $\{(a,c)\}$ for some $c \in S_3$ or $\{(a,c_2),(a,c_4)\}$. Let $v$ be
the vertex such that $vv_1v_2$ is an internal face in $G$. Let $G_1$ be the near-triangulation obtained
by deleting the edge $v_1v_2$ having $C_1 = v_1,v,v_2,\ldots,v_l$ as the boundary. Let $L(v) =
\{c_2,c_3,c_4\}$ and assign the label $L_2'$ to the edge $vv_2$, where $L_2' = \{(c_2,c),(c_4,c)\}$ if 
$L_2 = \{(a,c)\}$ and $L_2' = \{(c_3,c_2),(c_3,c_4),(c_2,c_4),(c_4,c_2)\}$ otherwise. The list
assignment to the vertices of $C_1$ satisfies the induction hypothesis, and by induction, there exists
a label $L_1$ that can be assigned to $v_1v_l$ such that $(L_1,L_2',L_3)$ is a feasible triple for
$G_1$. Since any $L$-colouring of $C_1$ is also an $L$-colouring of $C$, $(L_1,L_2,L_3)$ is a 
feasible triple for $G$. A symmetrical argument holds if $L(v_l) = \{b\}$ for some $b \in S_2$.

Suppose $i \in (1,l)$. If $v_i$ is of type $T_1$ then $L(v_i) \cap L(v_{i-1}) = \emptyset$, and if it is of 
type $T_2$, $L(v_i) \cap L(v_{i+1}) = \emptyset$, a contradiction.  If it is of type $T_3$ and
$v_{i-1}$ is also of type $T_3$, the same argument holds. The only possibility is that $i = p+1$, 
$L(v_{p+1}) = \{c_3\}$ and $L(v_{p}) = \{c_2,c_3,c_4\}$. The label $L_2$ assigned to the edge $v_pv_{p+1}$ 
must be $\{(c_2,c_3),(c_4,c_3)\}$. Let $v$ be the vertex such that $vv_pv_{p+1}$ is an internal face in $G$.
Let $G_1$ be the near-triangulation obtained by deleting the edge $v_pv_{p+1}$ with the cycle $C_1 =
v_1, \ldots,v_p,v,v_{p+1},\ldots,v_l$ as the boundary. Let $L(v) = \{c_1,c_2,c_4\}$ and assign the 
label $L_2' = \{(c_2,c_1),(c_4,c_1),(c_2,c_4),(c_4,c_2)\}$ to the edge $v_pv$. The list assignment to
the vertices of $C_1$, along with the label $L_2'$, satisfies the induction hypothesis, and by induction
there exists a feasible triple $(L_1,L_2',L_3)$ for $G_1$. In any $L$-colouring of $C_1$, the label
$L_2'$ ensures that $v_p$ is coloured $c_2$ or $c_4$, hence it is also an $L$-colouring of $C$.
Therefore $(L_1,L_2,L_3)$ is a feasible triple for $G$.

\noindent {\bf Case 2.2.2} We may now assume $|L(v_i)| = 3$ for all $i \in [1,l]$. 

\noindent {\bf Case 2.2.2.1} Suppose $p > 1$ and $q > p+1$. 
Suppose the label $L_2$ contains the pairs $(c_2,c_1)$ and $(c_4,c_1)$. Let $v_{p-1} = u_1,\ldots,u_r = 
v_{p+1}$ be the neighbours of $v_p$. Consider the near-triangulation $G_1$ obtained by deleting $v_p$ having
the cycle $v_1,\ldots,v_{p-1},u_2,\ldots,u_{r-1},v_{p+1}, \ldots, v_l$ as the boundary. Remove the elements 
$c_2,c_4$ from $L(v_{p+1})$ so that $L(v_{p+1}) = \{c_1\}$, which  converts $v_{p+1}$ to a vertex of type $T_1$. 
Let $L(u_i) = \{c_1,c_5\} \setminus L(u_{i+1})$ for $2 \le i < r$. Assign the label $L_2' = \{(c_1,c_2),
(c_1,c_4)\}$ to the edge $v_{p+1}v_{p+2}$ and keep the label $L_3$ for the edge $v_qv_{q+1}$. The resulting 
list assignment to the vertices of $C_1$ satisfies the induction hypothesis, and by induction, there exists a 
feasible triple $(L_1,L_2',L_3)$ for $G_1$. Since any $L$-colouring of $C_1$ can be extended to an 
$L$-colouring of $C$ by assigning $v_p$ a colour in $\{c_2,c_4\}$ that is not assigned to $v_{p-1}$, and since 
$(c_2,c_1),(c_4,c_1) \in L_2$, $(L_1,L_2,L_3)$ is a feasible triple for $G$.

If $(c_2,c_1) \not\in L_2$ then both $(c_3,c_2)$ and $(c_3,c_4)$ are in $L_2$. Now let $v_p = u_1,
\ldots,u_r = v_{p+2}$ be the neighbours of $v_{p+1}$. Let $G_1$ be the near-triangulation obtained by 
deleting the vertex $v_{p+1}$ with the cycle $C_1 = v_1,\ldots,v_p,u_2,\ldots,u_{r-1},v_{p+1},\ldots,v_l$
as the boundary. Remove the elements $c_2,c_4$ from $L(v_p)$ so that $L(v_p) = \{c_3\}$, which converts
$v_p$ to a vertex of type $T_3$. Let $L(u_i) = \{c_3,c_5\} \setminus L(u_{i-1})$ for $2 \le i < r$.
Assign the label $L_2' = \{(c_2,c_3),(c_4,c_3)\}$ to the edge $v_{p-1}v_{p}$ and retain the label $L_3$
for the edge $v_qv_{q+1}$ in $C_1$. The list assignment to the vertices of $C_1$ satisfies the induction
hypothesis and by induction there exists a feasible triple $(L_1,L_2',L_3)$ for $G_1$. Any $L$-colouring
of $C_1$ can be extended to an $L$-colouring of $C$ by assigning $v_{p+1}$ a colour in $\{c_2,c_4\}$
that is not assigned to $v_{p+2}$. This gives an $L$-colouring of $C$ in which $v_p$ is coloured
$c_3$ and $v_{p+1}$ is coloured $c_2$ or $c_4$, hence $(L_1,L_2,L_3)$ is a feasible triple for $G$.

\noindent {\bf Case 2.2.2.2} Suppose $p > 1$ and $q = p+1$. The argument is
similar to that in Case 2.2.2.1. Suppose $(c_2,c_1) \in L_2$, and again consider the near-triangulation obtained
by deleting the vertex $v_p$. If $(c_1,c_3) \in L_3$, assign $L(v_{p+1}) = \{c_1\}$ and $L(u_i) =
\{c_1,c_5\} \setminus L(u_{i+1})$ for $2 \le i < r$. There is no vertex of type $T_3$ in $C_1$ now. Instead of 
the label $L_3$ assigned to the edge $v_{p+1}v_{p+2}$, assign the corresponding label $L_3' = \{(c_1,c_3)\}$. 
Applying Lemma~\ref{two}, there exists a feasible pair $(L_1,L_3')$ for $G_1$. Assigning $v_p$ a colour
in $\{c_2,c_4\}$ that is not assigned to $v_{p-1}$ extends any $L$-colouring of $C_1$ to an $L$-colouring
of $C$. Therefore $(L_1,L_2,L_3)$ is a feasible triple for $G$. If $(c_1,c_3) \not\in L_3$, then $L_3 = 
\{(c_2,c_1),(c_2,c_4),(c_1,c_4),(c_4,c_1)\}$ and the same argument holds by choosing $L_3' = \{(c_1,c_4)\}$. 
The remaining possibility is $(c_2,c_1) \not\in L_2$, and by symmetry, $(c_2,c_1) \not\in L_3$ which implies 
$L_2 = \{(c_3,c_2),(c_3,c_4),(c_2,c_4),(c_4,c_2)\}$ and $L_3 = \{(c_1,c_3),(c_4,c_3),(c_1,c_4),(c_4,c_1)\}$. 
In this case, swap the colours $c_3$ and $c_4$ in all the lists and labels and assign list $\{c_3\}$ to the 
vertex $v_{p+1}$. The lists for the other vertices remain the same and the label $L_2$ changes to $\{(c_2,c_3),
(c_4,c_3)\}$ and $L_3$ to $\{(c_3,c_1),(c_3,c_4)\}$. This reduces to the Case 2.2.1 considered previously.

Cases 2.2.2.1 and 2.2.2.2 cover all possibilities with $p > 1$. A symmetrical argument can be used if 
$q < l-1$ by reversing the cycle and swapping the colours $c_1$ and $c_2$. 

\noindent {\bf Case 2.2.2.3} The only case that remains now is when $p=1$ and $q = l-1$, which means $v_1$ is 
the only vertex of type $T_1$ and $v_l$ is the only vertex of type $T_2$.

Suppose $l \ge 5$ and let $v_2 = u_1,u_2,\ldots,u_r = v_4$ be the neighbours of $v_3$ such that $v_3u_iu_{i+1}$ 
is an internal face in $G$, for $1 \le i < r$. Let $G_1$ be the near-triangulation obtained by deleting the 
vertex $v_3$ with the cycle $C_1 = v_1,v_2,u_2,\ldots,u_{r-1},v_4,\ldots,v_l$ as the boundary. Let $L(u_i) = 
\{c_3\}$ if $i$ is odd and $L(u_i) = \{c_5\}$ if $i$ is even for $2 \le i < r$. The list assignment to the 
vertices of $C_1$, together with the labels $L_2$ and $L_3$ for the edges $v_1v_2$ and $v_{l-1}v_l$ 
respectively, satisfies the induction hypothesis, and by induction there exists a feasible triple 
$(L_1,L_2,L_3)$ for $G_1$. Any $L$-colouring of $C_1$ can be extended to an $L$-colouring of $C$ by assigning 
$v_3$ a colour in $\{c_1,c_2,c_4\}$ that is not assigned to $v_2$ or $v_4$. This implies $(L_1,L_2,L_3)$ is 
also a feasible triple for $G$.

Suppose $l = 4$, $(c_2,c_1),(c_4,c_1) \in L_2$ and $(c_1,c_3),(c_4,c_3) \in L_3$. Let $v_4 = u_1,
u_2,\ldots,u_r = v_2$ be the neighbours of $v_1$ such that $v_1u_iu_{i+1}$ is an internal face in
$G$ for $1 \le i < r$.  Consider the near-triangulation $G_1$ obtained by deleting the vertex $v_1$ with
the cycle $C_1 = u_2,\ldots, v_2,v_3,v_4$ as the boundary. Assign the list $L(u_r) = L(v_2) = \{c_1\}$ 
and $L(u_i) = \{c_1,c_5\} \setminus L(u_{i+1})$ for $2 \le i < r$.  Assign the list $\{c_2,c_3,c_4\}$ to 
$v_3$ and $v_4$ so that all vertices in $C_1$ are of type $T_1$. By Lemma~\ref{one}, there exists an 
$L$-colouring $f$ of $C_1$ that extends to a $k$-colouring of $G_1$ with $f(v_4) = c_3$, $f(v_3) = c_4$
and $f(v_2) = c_1$. This can be extended to an $L$-colouring of $C$ by assigning the colour $c_2$ or $c_4$ 
to the vertex $v_1$. Then with $L_1 = \{(c_2,c_3),(c_4,c_3)\}$, $(L_1,L_2,L_3)$ is a feasible triple for $G$.
A symmetrical argument can be used if $\{(c_3,c_2),(c_3,c_4)\} \in L_2$ and $\{(c_2,c_1),(c_2,c_4)\} \in
L_3$.

Suppose $(c_2,c_1),(c_4,c_1) \in L_2$  but $(c_1,c_3),(c_4,c_3) \not\in L_3$. This implies 
$L_3 = \{(c_2,c_1), (c_2,c_4),(c_1,c_4),(c_4,c_1)\}$ and by the previous argument, $(c_3,c_2),(c_3,c_4)
\not\in L_2$. This implies $L_2 = \{(c_2,c_1),(c_4,c_1),(c_4,c_2),(c_2,c_4)\}$. In this case, the previous 
argument gives an $L$-colouring of $C_1$ that extends to a $k$-colouring of $G_1$ with $f(v_4) = c_4$, 
$f(v_3) = c_2$ and $f(v_1) = c_1$. Assigning colour $c_2$ to $v_1$ gives an $L$-colouring $f$ of $C$. Again 
using symmetry after reversing the cycle and swapping colours $c_1,c_2$, we get an $L$-colouring $f$ of $C$ 
that extends to a $k$-colouring of $G$ with $f(v_1) = c_4$, $f(v_2) = c_1$, $f(v_3) = c_2$ and $f(v_4) = c_1$. 
Therefore setting $L_1 = \{(c_2,c_4),(c_4,c_1)\}$ gives a feasible triple for $G$.

The remaining case to be considered is if $(c_2,c_1),(c_4,c_1) \not\in L_2$, and by symmetry,
$(c_2,c_1),(c_2,c_4) \not\in L_3$. Then $L_2 = \{(c_3,c_2),(c_3,c_4),(c_2,c_4),(c_4,c_2)\}$ and
$L_3 = \{(c_1,c_3),(c_4,c_3),(c_1,c_4),(c_4,c_1)\}$. Consider the  triangulation $G_1$,
and assign the lists $L(v_4) = \{c_1\}$, $L(u_i) = \{c_1,c_5\} \setminus L(u_{i-1})$ for
$2 \le i < r$, and  $L(v_2) =  L(v_3) = \{c_2,c_3,c_4\}$. Then all vertices are of type $T_1$ and there
exists an $L$-colouring $f$ of $C_1$ that extends to a $k$-colouring of $G_1$ with $f(v_4) = c_1$,
$f(v_2) = c_2$ and $f(v_3) = c_4$. Assigning $f(v_1) = c_3$ gives an $L$-colouring $f$ of $C$
that extends to a $k$-colouring of $G$.  A symmetrical argument, after reversing the cycle and
swapping the colours $c_1,c_2$ gives an $L$-colouring $f$ of $C$ that extends to a $k$-colouring 
of $G$ with $f(v_1) = c_2$, $f(v_2) = c_4$, $f(v_3) = c_1$ and $f(v_4) = c_3$. This implies
setting $L_1 = \{(c_3,c_1),(c_2,c_3)\}$ gives a feasible triple of labels $(L_1,L_2,L_3)$. This
completes all cases when $l = 4$ and hence the proof of Lemma~\ref{three}.
\hfill $\Box$

\noindent {\bf Proof of Theorem~\ref{main}: } The proof follows from Lemmas~\ref{one}, \ref{two} and
\ref{three}. Let $f$ be a good $k$-colouring of a cycle $C$ of length $l$. Suppose $f$ satisfies condition 1 
in the definition of a good colouring. By permuting the colours, we can assume that $F[1,l] \subseteq S_1$.
Then assigning the list $L(v_i) = \{f(v_i)\}$ for all $1 \le i \le l$ gives a list assignment to the 
vertices of $C$ that satisfies the conditions of Lemma~\ref{one}. Since any chordless near-triangulation
$G$ with boundary $C$ is consistent with the list assignment to $C$, Lemma~\ref{one} implies the colouring of 
$C$ can be extended to a $k$-colouring of $G$.

Suppose $f$ satisfies condition 2 in the definition of a good colouring, with indices $p$ and $q$.
By relabeling the vertices and permuting the colours, we may assume that $p=1$, $F[1,q) \subseteq 
S_1$, $F[q,1) \subseteq S_2$, and the set $A = \{c_j | 5 \le j \le k\}$. Then the list
assignment $L(v_i) = \{f(v_i)\}$ for $1 \le i \le l$, with the label $L_2 = L(v_{q-1}) \times L(v_q)$
assigned to the edge $v_{q-1}v_q$, satisfies the conditions of Lemma~\ref{two} and the result follows.

Suppose $f$ satisfies condition 3 in the definition of a good colouring, with indices $p, q, r$.
Again by relabeling vertices and permuting colours, we may assume that $p = 1$, $F[1,q)
\subseteq S_1$, $F[q,r) \subseteq S_3$ and $F[r,1) \subseteq S_2$, where again $A = \{c_j | 5 \le j \le k\}$. 
The list assignment $L(v_i) = \{f(v_i)\}$  satisfies the conditions of Lemma~\ref{three} with the labels 
$L_2 = L(v_{q-1}) \times L(v_q)$ and $L_3 = L(v_{r-1}) \times L(v_r)$ assigned to the edges $v_{q-1}v_q$
and $v_{r-1}v_r$, respectively. The result follows from Lemma~\ref{three}.

This completes the proof of Theorem~\ref{main}.
\hfill $\Box$  

\begin{cor}
\label{longest}
Any $k$-colouring $f$ of a cycle of length $l \le 2k-5$ such that $|F[1,l]| \le
k-1$ is safe. There exists a $k$-colouring of a cycle of length $2k-4$ that
is not safe and uses $k-1$ distinct colours.
\end{cor}

\noindent {\bf Proof:} If $|F[1,l]| \le k-3$ then $f$ is good and the result
follows from Theorem~\ref{main}. If $|F[1,l]| = k-2$, we may assume $F[1,l]
= \{c_1,\ldots,c_{k-2}\}$. Since $l \le 2k-5$, at least one colour occurs
only once in the cycle. Without loss of generality, by relabeling vertices
and permuting colours, we may assume $v_l$ is the only vertex with colour
$c_{k-2}$. Then $f$ satisfies condition 2 in the definition of a good
$k$-colouring, with $p = 1$, $q =l$ and $A = \{c_1,\ldots,c_{k-4}\}$. If $|F[1,l]| = k-1$, 
there are at least three distinct colours that occur exactly once in the colouring of the
cycle. Without loss of generality, we may assume these colours are $c_{k-3},c_{k-2},c_{k-1}$ and 
the colour $c_k$ does not occur. Let $1 < p' < q' \le l$ be indices such that $f(v_1) = c_{k-3}$, 
$f(v_{p'}) = c_{k-2}$ and $f(v_{q'}) = c_{k-1}$. Then $f$ satisfies condition 3 in the definition of a 
good colouring with $p = 1, q = p', r = q'$ and $A = \{c_1,\ldots,c_{k-4}\}$. Theorem~\ref{main} implies $f$ 
is safe.

An example of a $k$-colouring $f$ of a cycle of length $2k-4$ that is not
safe is given by $f(v_i) = f(v_{k-2+i}) = c_i$ for $1 \le i \le k-3$,
$f(v_{k-2}) = c_{k-2}$ and $f(v_{2k-4}) = c_{k-1}$. Note that this
colouring satisfies condition 2 in Lemma~\ref{unsafe}, and is bad.
\hfill $\Box$

If Conjecture~\ref{safe} is true, then every $k$-colouring of a cycle
of length at most $3k-10$ that uses at most $k-2$ distinct colours is
safe. However, Theorem~\ref{main} does not prove this since the
$k$-colouring of a cycle of length $2k-4$  defined by $f(v_i) = f(v_{i+k-2}) =
c_i$ for $1 \le i \le k-2$ is not good.

\section {Remarks}

In principle, it may be possible to extend this approach to prove the
four colour theorem itself. However, it is not sufficient to use 
a set of $4$-colourings of a cycle defined by list-colourings. A
more general way of defining sets of colourings is required. A set 
$\mathcal{F}$ of $4$-colourings of a cycle $C$ of length $l$ is safe, if for 
every chordless near-triangulation $G$ with $C$ as the boundary, there exists
a $4$-colouring $f \in \mathcal{F}$ of $C$ that can be extended to a $4$-colouring
of $G$. While it would be nice to characterize exactly the safe sets of
$4$-colorings, proving the four colour theorem only requires showing
that any proper colouring of a triangle is safe. It may be possible to
do this by finding a simpler sufficient condition that ensures safety of
a set of colourings.

\end{document}